\documentclass[a4paper,10pt]{article}
\usepackage{amssymb}
\usepackage{amsmath}
\usepackage{amsthm}
\usepackage{graphicx}

\newtheorem{thm}{Theorem}[section]

\newtheorem{lem}[thm]{Lemma}
\newtheorem{prop}[thm]{Proposition}

\theoremstyle{definition}
\newtheorem{defn}[thm]{Definition}
\newtheorem{rem}[thm]{Remark}
\newtheorem{ej}[thm]{Example}

\def \Cram{{\rm Cram}}

\title{Tropical constructive Pappus' theorem\footnote{
First published in Int. Math. Res. Not, 2005(2005) nº 39, published by Hindawi Publishing Corporation}}
\author{Luis Felipe Tabera
\footnote{Luis Felipe Tabera (luisfelipe.tabera@unican.es) is
supported by the European Research Training Network RAAG
(HPRN-CT-2001-00271) and also by a Formaci\'on de
Profesorado Universitario research grant from the Spanish
Ministerio de Educaci\'on y Ciencia.}
}

\begin{document}
\maketitle
\begin{abstract}
In this paper, we state a correspondence between classical
and tropical Cramer's rule. This correspondence will allow
us to compare linear geometric constructions in the
projective and tropical spaces. In particular, we prove a
constructive version of Pappus' theorem, as conjectured in
\cite{RGST}.
\end{abstract}

\section{Introduction}
In the last years we have seen an increasing interest in
tropical geometry. Introductory papers in tropical geometry
may be found by the interested reader in Richter-Gebert et
al. \cite{RGST}, the specific chapter in the book of
Sturmfels \cite{STU} or the survey due to Mikhalkin
\cite{MIK2}. In the last reference \cite{MIK2}, G. Mikhalkin
applies tropical geometry to enumerative geometry, proving a
new way to calculate Gromov-Witten invariants in the
projective plane. These invariants can be used to count the
number of curves, with given genus and degree, passing through
a configuration of points. This method was first suggested
by Kontsevitch and it is also approached from another point
of view in \cite{SHU}. Also, in \cite{ITKASHU}, we can find
some computations of bounds for the Welschinger invariant in
several toric surfaces using tropical geometry, which are
interpreted as the algebraic count of real rational curves
through a real configuration of points. Moreover, we can see
an application of tropical geometry to combinatorics in
\cite{SPESTU}. Thus we observe that tropical geometry is a
powerful tool to study different branches of mathematics. The
problem is that it is not easy to translate familiar
geometric definitions to a tropical framework.

This paper deals with the specific problem of successfully
translating Pappus theorem to a suitable tropical state,
using the notion of stable intersection and stable join, as
presented in \cite{RGST}. A more systematic study of tropical
constructions and their relationship with classical ones is
treated in \cite{Tab}.

We will work on the tropical semiring
$(\mathbb{T},\oplus,\odot) = (\mathbb{R}, \max, +)$, the set
of real numbers with the tropical addition $a \oplus b =
\max\{a,b\}$ and the tropical product $a \odot b = a + b$.
In \cite{RGST} and other references, it is also used the
tropical semiring $(\mathbb{R}, \min, +)$ instead. But it is
straightforward to check that all the results can be
translated from one point of view to the other using the
isomorphism $a\mapsto -a$.

Now, we present the objects that we will work with. Given a
tropical polynomial $f = \bigoplus_{i\in I}a_i\odot x^i$,
where $i=(i_1,\ldots,i_n)$ and $x^i=x_1^{i_1}\odot \ldots
\odot x_n^{i_n}$, we define the \emph{ tropical variety
associated with} $f$ as the set $\mathcal{T}(f) :=
\{x=(x_1,\ldots,x_n)\in \mathbb{T}^n\ |\ f(x) = \max\{
a_i+x_1i_1 + \ldots + x_ni_n,\ i\in I\}$ \emph{is attained
for at least two different} $i\}$. That is, the set of
points where $f$ is not differentiable.

In the following, we will use homogeneous coordinates in the
tropical space $\mathbb{T}^n$, representing the point
$(y_1,\ldots,y_n)\in \mathbb{T}^n$ by $[y_1:\ldots:y_n:0]$,
with the identification $[y_1:\ldots:y_{n+1}]=[\alpha\odot
y_1:\ldots:\alpha\odot y_{n+1}]=[ \alpha +y_1:\ldots :\alpha+
y_{n+1}]$, $\alpha\in\mathbb{T}$.
We recover the affine coordinates using the usual
subtraction (there is no notion of tropical subtraction),
$[y_1 : \ldots : y_n : y_{n+1}] = (y_1-y_{n+1}, \ldots,
y_n-y_{n+1})$. We use homogeneous coordinates because it is
easier to state Cramer's rule in this context.

In this direction, let us remark that the simplest well
known varieties are tropical lines in the plane.
For instance, take
$f = a\odot x\oplus b\odot y\oplus c\odot z$ a linear
homogeneous polynomial. The corresponding tropical line in
$\mathbb{T}^2$ is the set $[x:y:z]$ such that $(a+x=b+y\geq
c+z)$ or $(a+x=c+z\geq b+y)$ or $(b+y=c+z\geq a+x)$.
We obtain three rays emerging from the point $[-a:-b:-c]$,
with vectors in the directions $[0:0:-1]$, $[0:-1:0]$,
$[-1:0:0]$.

\begin{figure}
\begin{center}
\begin{picture}(70,80)
\put(0,40){\line(1,0){50}}
\put(25,40){\line(0,-1){40}}
\put(25,40){\line(1,1){40}}
\put(50,40){\line(0,-1){40}}
\put(50,40){\line(1,1){40}}
\linethickness{1.2pt}
\put(0,40){\line(1,0){25}}
\end{picture}
\caption{An infinite intersection of two lines}\label{interinfinita}
\end{center}
\end{figure}
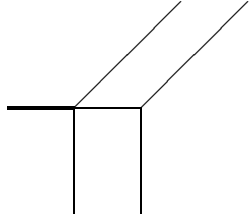

Now we arrive to the following problem. What should be
considered as the intersection of two given lines or, more
generally, the intersection of tropical hypersurfaces. It is
not trivial, as it may happen that two different lines share
an infinite number of points, see figure
(\ref{interinfinita}).
An answer is the following: given two tropical lines, there
exists only one point in the intersection such that it is
stable (in some sense) under small perturbations of the two
lines, see \cite{RGST}. This distinguished point is called
the \emph{stable intersection} of the lines. Similarly, given
two points, there is only one line that passes through the
two given points and is stable under small perturbation of
the two points, it is called \emph{the stable line} or the
\emph{stable join} of the points. Both stable intersection
and stable join can be computed using the tropical analog of
Cramer's rule, as follows.

First, the tropical determinant of a given $n\times n$ matrix
in $\mathbb{T}$ is defined as
\[\begin{vmatrix}a_{11}&\cdots&a_{1n}\\
\vdots&\cdots&\vdots\\
a_{n1}&\cdots&a_{nn}\end{vmatrix}_t=\bigoplus_{\sigma
\in\Sigma_n} a_{1\sigma(1)} \odot\ldots\odot a_{n\sigma(n)}.\]
Suppose now a tropical linear system of $n$ equations in
$n+1$ homogeneous variables is given. We write $O$ the
$n\times (n\!+\!1)$ matrix of coefficients and denote by
$O^{i}$ the matrix resulting from deleting the $i$-th column
of $O$. Then, it is shown in \cite{RGST} that the point
$[ | O^{1} |_t : \ldots : |O^{n\!+\!1} |_t ]$
is not only a common point of the $n$ hyperplanes, but the
only one which is stable under small perturbations in the
coefficients of the hyperplanes. Therefore, this version of
Cramer's rule is the right tool that gives us the stable
intersection of $n$ hyperplanes in the tropical space
$\mathbb{T}^n$. The reader might guess that it is also the
tool to compute the stable join of $n$ points, as this
problem can also be seen as solving a regular linear system
of equations, the unknowns are the coefficients $a_i$ of the
linear equation defining the hyperplane
$\bigoplus_{i=1}^{n+1}a_i\odot x_i$ and the row entries of
the matrix that represent the system are the (homogeneous)
coordinates of the given $n$ points.

Another definition of tropical hypersurfaces may be taken as
starting point following the idea that tropical varieties are
non-archimedean amoebas. After Gelfand, Kapranov and
Zelevinski \cite{GKZ}, an amoeba is the image by the
logarithmic function of an algebraic variety in
$(\mathbb{C}^*)^n$. As it is usual in tropical geometry, we
will work instead with the field $\mathbb{K}$ of ``Puiseux''
series with complex coordinates and real exponents. Its
elements are formal series $\sum_{i\in \Lambda}\alpha_it^i$
where $\alpha_i\in\mathbb{C}$ and $\Lambda\subset\mathbb{R}$
is a countable set contained in a finite number of arithmetic
sequences. Moreover we use, over $\mathbb{K}$, its
non-archimedian valuation in order to define our amoebas.

For this purpose we define in $\mathbb{K}$ the application
$T :\mathbb{K}^* \longrightarrow \mathbb{R}$, $T(x)=-o(x)$,
(i.e. minus the order of the Puiseux series). What we obtain
is, in fact, an application
$T:\mathbb{K}^*\longrightarrow\mathbb{T}$ onto our tropical
semiring. This application satisfies that $T(xy)=T(x)\odot
T(y)$ and $T(x+y)=T(x)\oplus T(y)$ if $T(x)\neq T(y)$.

Now, let $V$ be an algebraic variety in the algebraic torus
$(\mathbb{K}^*)^n$. The result of applying $T$ component-wise
over $V$, $T(V)$ is called, by definition, \emph{an algebraic
tropical variety} and also the \emph{tropicalization} of $V$.
Conversely, for an algebraic tropical variety $U$, a
\emph{lift of $U$} is any algebraic variety $V$ such that
$T(V)=U$.

Let us remark that we have given two different definitions
of tropical hypersurfaces: First, as the variety
$\mathcal{T}(f)$ associated to a tropical polynomial $f$;
second, as the tropicalization $T(H)$ of an algebraic
hypersurface $H$. The following theorem of Kapranov
\cite{EKL} shows that the two different points of view yield,
for hypersurfaces, the same subsets. Precisely, let
$\widetilde{f}=\sum_{i\in I}a_ix^i$ be a polynomial in
$\mathbb{K}[x_1,\ldots,x_n]$, $a_i\neq 0$, $i\in I$. Let
$f=\bigoplus_{i\in I}T(a_i)\odot x^i$ be the corresponding
tropical polynomial. Then, $\mathcal{T}(f)$ equals
$T(\{\widetilde{f}=0\})$.

Notice that this theorem does not hold for general
varieties, not even for complete intersections. It is easy to
find examples of varieties such that $T(\{ \widetilde{f}_1=0
\})\cap\dots\cap V(\{\widetilde{f}_n=0\})) \neq
\mathcal{T}(f_1) \cap\ldots\cap\mathcal{T}(f_n)$ (for some
examples, cf. \cite{Tab}). Therefore, in order to study this
situation with more detail one possible way is trying to
enlarge the set of defining equations for the algebraic
variety so that its tropicalization coincides with the
intersection of the tropicalization of the given equations
(see \cite{RGST} for the particular case of linear
varieties). Another possibility, that we will follow here, is
to restrict the intersection of tropical hypersurfaces to the
particularly meaningful subset of stable points.

Correspondingly, a classical theorem of elementary geometry
can be regarded in at least two different ways: first,
describing by algebraic equations its hypotheses and thesis.
Second, as a collection of construction steps with
ge\-o\-me\-tric entities. This duality appears in
\cite{RGST}, where two versions of Pappus' theorem are
presented, showing that the tropicalization of the algebraic
translation of the hypotheses does not yield the thesis
without adding extra polynomial equations (thus obtaining a
tropical basis of the hypotheses ideal) in order to have a
correct translation to the tropical case.

On the other hand, \cite{RGST} formulates a conjecture about
the validity of the straightforward translation to the
tropical context of a constructive version of Pappus'
theorem. Here, of course, one must keep up, without
modification in the tropical framework, with the given
collection of construction steps. The goal of this paper is
precisely to prove this conjecture. Let us introduce some
notation needed to state it.

First of all, using duality, we identify the line $a\odot
x\oplus b\odot y\oplus c\odot z$ with the point $[a:b:c]$ in
$\mathbb{T}^2$. The origin of the rays of the line defined by
the polynomial is the point $[-a:-b:-c]$. Now, if we have two
lines $[a:b:c]$, $[d:e:f]$, then the stable intersection
corresponds to the point:
$\left[\begin{vmatrix}b&c\\e&f\end{vmatrix}_t :
\begin{vmatrix}a&c\\d&f\end{vmatrix}_t :
\begin{vmatrix}a&b\\d&e\end{vmatrix}_t
\right]$, which is the stable solution of
the system $a\odot x\oplus b\odot y\oplus c\odot z$, $d\odot
x\oplus e\odot y\oplus f\odot z$. Also, if we have two
points $[a:b:c]$, $[d:e:f]$, the previous expression
corresponds to the coordinates of the stable line defined by
those points. Thus, as in \cite{RGST}, we define the cross
product of two points $x=[x_1:x_2:x_3]$ and $y=[y_1:y_2:y_3]$
as \[x\otimes y=\left[\begin{vmatrix}x_2&x_3 \\
y_2&y_3\end{vmatrix}_t : \begin{vmatrix}
x_1&x_3\\y_1&y_3\end{vmatrix}_t:
\begin{vmatrix}x_1&x_2\\y_1&y_2\end{vmatrix}_t\right ],\]
which can be interpreted as the intersection of two lines or
as finding the line through two points, depending on the
context.

With this terminology theorem \ref{Pappus} states that there exists a tropical construction such that given five points
$1$, $2$, $3$, $4$, $5$ in the tropical plane, it computes three additional points $6$, $7$, $8$ and nine lines $a$, $b$, $c$,
$a'$, $b'$, $c'$, $a''$, $b''$, $c''$ such that the hole set of elements is always in Pappus position (in the sense of \cite{RGST})
and hence, the intersection of $a''$, $b''$, $c''$ is not empty.

Remark that in the thesis of the theorem, we do not mean that
the three lines $a''$, $b''$, $c''$ share a point that is
stable under perturbations but just that their intersection
is non empty.

In our approach to proving this theorem it is essential to
understand the behavior of Cramer's rule. In \cite{RGST},
Cramer's rule is analyzed using generic small perturbations
in the coefficients of the system. Our point of view is a
little bit different, as we try to compare the performance of
Cramer's rule under the valuation map $T$. This allows us to
give sufficient conditions for a chain of computations of
tropical determinants to be lifted to the Puiseux series
field. Using this lift, we are able to derive results from
the classical to the tropical context.

The main idea is to take the input elements for the geometric
construction in the tropical space, make a lift to the
projective space over the Puiseux series field and perform
there the given construction using classical Cramer's rule.
Then we prove that, in this particular construction, if the
elements in the lift are taken general enough, the results
given by the application of Cramer's rule in the projective
ambient should correspond with those obtained applying
Cramer's rule in the tropical case.

Unfortunately, this procedure does not hold for other
constructions, see counterexample \ref{ejemplociclo}.
In fact, we can ensure that the good behavior in Pappus' case
happens because our construction is of a very particular
kind. Namely, we will prove that tropical and projective
constructions behave well with respect to tropicalization
when a certain graph associated with the construction is a
tree.

The paper is structured as follows. In section 2 we will
study with detail the relation between classical and tropical
constructions, including a new proof (\ref{pseudocramerlema})
of Cramer's rule. Then, we associate a graph to a tropical
construction (\ref{grafo}) and we introduce the notion of a
\emph{tropically admissible} construction of an elements
(\ref{tropadm}). Finally we state (\ref{Levantamientogeneral})
the validity of the specialization of a chain of Cramer's
rule computations to the tropical space, when the
construction graph is a tree. Section 3 is devoted
to prove the conjectured version of Pappus' theorem
(\ref{Pappus}), including some comments and remarks. We
conclude (Section 4) with some reflections on the difficulty
of achieving more general results on this topic.

\section{Tropical Geometric Constructions} By a
\emph{geometric construction in the classical case} we will
understand an abstract procedure consisting of
\begin{itemize}
\item Input data: A finite number of points or lines in
$\mathbb{P}^2(\mathbb{K})$ that will eventually specialize
to concrete elements given by its homogeneous coordinates (in
the case of lines, the coordinates of the corresponding
point in the dual plane).
\item Allowed steps: computing the
\begin{itemize}
\item line passing through two points
\item intersection point of two lines.
\end{itemize}
\item Output: A finite set of points and lines
\end{itemize}
Likewise a \emph{geometric construction in the tropical
plane} consists of a similar procedure, replacing in the
steps above the ``line through two points'' by the ``stable
line passing through two points'' and the ``intersection of
two lines'' by the ``stable intersection of two lines''.

We want to study the relation between a given construction in
the classical setting and the corresponding tropical one, see
\cite{Tab} for a more general study of tropical geometric
constructions. Namely, we want to analyze, for different
constructions, the commutativity of the following diagram:
\begin{equation}\label{diagrama}
\begin{matrix} (\mathbb{K}^*)^2&&\mathbb{T}^2\\ {\rm
Input}&\stackrel{T^{-1}}{\longleftarrow}&{\rm Input}\\
\downarrow&&\downarrow\\ {\rm Output}&
\stackrel{T}{\longrightarrow} &{\rm Output}
\end{matrix}
\end{equation}
where $T$ stands for the tropicalization mapping. That is,
given a construction, we want to study when, for some given
tropical input data, we are able to find a suitable lift of
the input data to the Puiseux series field $\mathbb{K}$,
perform the construction in that projective plane,
tropicalize all the output elements and find out that they
are exactly the elements obtained by the tropical
construction.

We will soon notice that it is not always possible. Even if
it holds for some constructions, it will not do for every
choice of an input lift (example \ref{restricciones}). Let us
start with the simplest case of one step constructions
involving Cramer's rule only once.

For a Puiseux series $S=\alpha t^{k}+\ldots$, $\alpha\neq 0$,
we will denote by $Pc(S)=\alpha$ the \emph{principal
coefficient} of the series.

Let $B=(b_{i,j})$ be a $n\times n$ matrix in $\mathbb{K}^*$.
Let us start by studying conditions for the commutativity
of tropicalization and determinant computation, i.e.
establishing when computing the determinant of $B$ and
tropicalizing it equals the determinant of $T(B)$. First, it
can be easily checked (\cite{Tab}) that this this equality
does not hold in general. Now, we will show that the
conditions to have this property can be expressed in terms
of the principal coefficients of the entries of $B$.

In this context we need to introduce the following
terminology.

\begin{defn}\label{pseudodeterminante} Let $O=(o_{ij})$ be a
$n\times n$ matrix with coefficients in $\mathbb{T}$. Let
$A=(a_{ij})$ be a $n \times n$ matrix in a ring $R$. We
denote by $|O|_t$ the tropical determinant of $O$ and we
define \[\Delta_O(A)=\sum_{\substack{\sigma\in\Sigma_n\\
o_{1,\sigma(1)}\odot \ldots\odot o_{n,\sigma(n)} = |O|_t}}
(-1)^{i(\sigma)}a_{1,\sigma(1)}\cdot\ldots\cdot
a_{n,\sigma(n)}\]
the \emph{pseudo-determinant of $A$} with respect to $O$.
\end{defn}

\begin{lem}\label{pseudodeterminantelema} Let $B=(b_{i,j})$
be a matrix in $\mathbb{K}^*$, $A=(a_{i,j})$ the matrix of
principal coefficients in $B$ and $O$ the tropicalization
matrix of $B$, $o_{i,j}=T(b_{i,j})$. If $\Delta_O(A)\neq 0$,
then the principal coefficient of $|B|$, $Pc(|B|)$ equals
$\Delta_O(A)$. Moreover $T(|B|)$ coincides with the tropical
determinant $|O|_t$.
\end{lem}
\begin{proof} Notice that in the expansion of $|B|$, the
permutations $\sigma$, where the order of the corresponding
summand $b_{1,\sigma(1)}\cdot\ldots\cdot b_{n,\sigma(n)}$
is the smallest possible one, are exactly the permutations
in the expansion of $|O|_t$ where $|O|_t$ is attained. So,
the coefficient of the term $t^{-|O|_t}$ is $\Delta_O(A)$. If
it is non zero, then the order of $|B|$ is $-|O|_t$.
\end{proof}

Now, we can extend this lemma to the context of Cramer's rule.

\begin{defn}\label{pseudocramer} Let $O=(o_{ij})$ be a
$n\times(n\!+\!1)$ tropical matrix. Let $A=(a_{ij})$ be a
matrix in a ring $R$ with the same dimension as $O$. We
define \[\Cram_O(A)=(S_1,\ldots,S_{n+1})\] where
$S_i=\Delta_{O^{i}}(A^{i})$ and $O^{i}$ (respectively,
$A^{i}$) denote the corresponding submatrices obtained by
deleting the $i$-th column in $O$ (respectively, $A$).
\end{defn}

\begin{lem}\label{pseudocramerlema} Suppose we are given a
linear equation system in the semiring $\mathbb{T}$, with
$n$ equations in $n\!+\!1$ homogeneous variables. Let $O$
be the coefficient matrix of the system. Let $B$ be any
matrix such that $T(B)=O$. Let $A$ be the principal
coefficient matrix of $B$. If no element of $\Cram_O(A)$
vanishes, then the linear system defined by $B$ has only
one projective solution and its tropicalization equals the
stable solution $[ | O^{1} |_t : \ldots : |O^{n\!+1\!} |_t ]$
\end{lem}
\begin{proof} Apply the previous lemma to every component
of the projective solution.
\end{proof}

\begin{prop}\label{consunpaso}
If we have a one step construction, namely the stable join of
two points or the stable intersection of two lines, then for
every specialization of the input data, there exists a
concrete lift that makes diagram (\ref{diagrama}) commutative.
\end{prop}
\begin{proof}
Lemma \ref{pseudocramerlema} gives us sufficient conditions
in the principal coefficients of the input elements lift to
assure that the diagram is commutative. As the
pseudodeterminants are nonzero polynomials, there is always a
possible choice of coordinates that makes all
pseudodeterminants nonzero.
\end{proof}

\begin{rem}
Lemma \ref{pseudocramerlema} is not only useful to compute
the intersection of $n$ hyperplanes in $\mathbb{T}^n$, but it
is also valid, for example, to compute the stable plane conic
through five points, as it can be also interpreted as finding
the stable intersection of 5 hyperplanes in the space of
tropical plane conics $\mathbb{T}^5$.
\end{rem}

\begin{ej}\label{restricciones} We take here the two lines
$2\odot x\oplus (-3)\odot y\oplus 0\odot z$ and
$(-4)\odot x\oplus (-3)\odot y\oplus 0\odot z$. Their
intersection is the set $\{[t:3:0]\ |\ t\leq -2\}$ which is
not the tropicalization of any variety. If we perform, as
above, tropical Cramer's rule on these data, we obtain the
point $[-2:3:0]$. This is the limit of the unique
intersection point under small generic perturbations of the
tropical lines. If we take lifts of these lines,
they take the form $\alpha_xt^{-2}x + \alpha_yt^3y +
\alpha_zz$, $\beta_xt^4x + \beta_yt^3y + \beta_zz$, where we
do not write higher order terms. The solution of this
system is $[(\alpha_y\beta_z - \beta_y\alpha_z)t^3 :
-(\alpha_x\beta_zt^{-2} -\alpha_z\beta_xt^4) :
\alpha_x\beta_yt-\alpha_y\beta_xt^7]$.

We apply the notation of \ref{pseudodeterminante} over the
first coordinate data, yielding:\\
$A=\begin{pmatrix} \alpha_y &\alpha_z \\
\beta_y&\beta_z\end{pmatrix}$,
$O=\begin{pmatrix}-3&0\\-3&0\end{pmatrix}$, $\Delta_O(A) =
\alpha_y\beta_z-\beta_y\alpha_z$. Thus, if $\Delta_O(A)\neq
0$, then the intersection point will tropicalize to
$[-3:2:-1]=[-2:3:0]$. But if $\Delta_O(A) = 0$, we cannot know
what the tropicalization of the intersection is, since it
depends on the series higher order terms that we have not
written. We only know that they will be of the form $[t:3:0]$,
$t\leq -2$. It can even be the case that the intersection
has its first coordinate equals to 0. In this case there
will be no tropicalization of the intersection at all.
\end{ej}

Therefore we cannot expect, in general, that the
tropicalization of the generators of an ideal will describe
the tropicalization of the variety this ideal generates.

Now we proceed with the case of more complicated
constructions. Namely, for those constructions such that some
elements are computed from previously constructed elements.
In this respect the following notation and lemmas are useful.

\begin{defn}\label{grafo}
Given a linear geometric construction, we associate a graph
to every element $P$ of the construction, that is, to all the
points and lines appearing at some step of the construction,
including input, intermediate and output elements. The
vertices of the graph associated to an element $P$ will
correspond to all the elements that we have recursively used
to construct $P$. We link every element with the elements
from which it is constructed directly. That is, if point $a$
(respectively line $a$) is the intersection of points $b$,
$c$ (respectively the join of points $b$ and $c$), then we
write edges $ab$ and $ac$. We call this graph the
\emph{construction graph of $P$}.

For instance, the construction graph of and input element
consists in just one point representing the element. Other
examples appear in \ref{Pappus} or in \ref{ejemplociclo}.
\end{defn}

\begin{defn}\label{tropadm}
We say that the construction of an element $P$ is
\emph{tropically admissible (by Cramer's rule)} if its
associated construction graph is a tree.
\end{defn}

In general, given a tropical construction and some concrete
input data, we will study the possible lifts of these input
elements to the projective setting, by parametrizing the
principal coefficients of the corresponding series with
different variables. The following lemma describes the effect
of the corresponding classical construction on this generic
lift.

\begin{lem}\label{Levantamientounpaso}
Let $C_i=\{c_i^1,\ldots,c_i^{j_i}\}$, $1 \leq i\leq k$ be
disjoint sets of variables.
Suppose that we have $F_u = \{f_u^1,\ldots,f_u^{n+1}\}
\subseteq \mathbb{C}[\bigcup_{i=1}^k C_i]$, $1\leq u \leq n$
sets of polynomials in the variables $c_i^j$.
Suppose also that the following properties hold:
\begin{itemize}
\item For a fixed set $F_u$, $f_u^l$, with $1\leq l\leq n+1$
are multihomogeneous polynomials in the sets of variables
$C_{u^1},\ldots,C_{u^{s_u}}$ with the same multidegree.
\item If $u\neq v$ then $F_u$, $F_v$ involve different sets
of variables $C_i$.
\item In a family $F_u$, if $l\neq m$ then the monomials of
$f_u^l$ are all different from the monomials of $f_u^m$.
\end{itemize}
Let us construct the $n\times (n\!+\!1)$ matrix
\[A=(f_u^l)_{\substack{1\leq u\leq n,\\ 1\leq l\leq n+1}}\]
and suppose that we are given a $n\times (n\!+\!1)$ matrix
$O$ in $\mathbb{T}$.
Write \[S=\Cram_O(A)=(S_1,\ldots,S_{n+1}).\] Then
\begin{enumerate}
\item $S_1,\ldots,S_{n+1}$ are non-zero multihomogeneous
polynomials in the sets of variables $C_1,\ldots,C_k$ with
the same multidegree.
\item If $\sigma, \tau$ are different permutations in
$\Sigma_{n+1}$ which appear in the expansion of $S_l$ (and,
therefore $\sigma(n+1)=\tau(n+1)=l$), then all resulting
monomials in $\prod_{u=1}^n (A^{l})_u^{\sigma(u)}$ are
different from the monomials in
$\prod_{u=1}^n (A^{l})_u^{\tau(u)}$
\item If $l\neq m$, then $S_l$, $S_m$ have no common
monomials.
\end{enumerate}
\end{lem}
\begin{proof}
First we prove {\emph 2}. If we have two different
permutations $\sigma$, $\tau$, there is a natural number $v$,
$1\leq v\leq n$ where the permutations differ, then the
monomials in $f_v^{\sigma(v)}$, $f_v^{\tau(v)}$ are all
different and these polynomials are the only factors of the
products $\prod_{u=1}^n (A^{l})_u^{\sigma(u)}$,
$\prod_{u=1}^n (A^{l})_u^{\tau(u)}$ where we find the
variables which appear in the family $F_v$. It follows that
these products cannot share any monomial. In particular, in
the sum of several of these products, there is no
cancellation of monomials, proving item {\emph 1}. So, in
fact, we obtain that different minors share no monomial
and we obtain immediately {\emph 3}. All those minors must
have the same multidegree, which is just the concatenation of
the multidegree of the family $F_1,\ldots,F_n$, by
construction.
\end{proof}

\begin{ej}
At this point it may be helpful to give an example of the
lemma.
Consider the sets\\
$C_1=\{x,y\}$, $C_2=\{z\}$, $C_3=\{m,n\}$, $C_4=\{o,p,q\}$, $C_5=\{r\}$.\\
$F_1=\{x^2yz+y^3z,x^3z,2xy^2z\}$\\
$F_2=\{mnor^2,m^2or^2+mnpr^2,n^2or^2+m^2pr^2+n^2pr^2\}$\\
Every polynomial in $F_1$ is multihomogeneous in $C_1$, $C_2$
with multidegree $(3,1)$.\\
Every polynomial in $F_2$ is multihomogeneous in $C_3$,
$C_4$, $C_5$ with multidegree $(2,1,2)$.\\
All the monomials in the polynomial are different.\\
Then, matrix $A=\begin{pmatrix} x^2yz+y^3z&x^3z&2xy^2z\\
mnor^2&m^2or^2+mnpr^2&n^2or^2+m^2pr^2+n^2pr^2\end{pmatrix}$.\\
We take as matrix $O$ in $\Cram_O(A)$,
$O=\begin{pmatrix}1&2&3\\0&3&2\end{pmatrix}$\\
$S_1=(m^2or^2+mnpr^2)(2xy^2z)=2xy^2zm^2or^2+2xy^2zmnpr^2$\\
$S_2=(x^2yz+y^3z)(n^2or^2+m^2pr^2+n^2pr^2)+(mnor^2)(2xy^2z)=
x^2yzn^2or^2+x^2yzm^2pr^2+x^2yzn^2pr^2+y^3zn^2or^2+y^3zm^2pr^2
+y^3zn^2pr^2+2xy^2zmnor^2$\\
$S_3=(x^2yz+y^3z)(m^2or^2+mnpr^2)=x^2yzm^2or^2+x^2yzmnpr^2+
y^3zm^2or^2+y^3zmnpr^2$.\\
Finally, we check that the polynomials $S_1$, $S_2$, $S_3$
share no monomial and are multihomogeneous in $C_1$, $C_2$,
$C_3$, $C_4$, $C_5$ with multidegree $(3,1,2,1,2)$.
\end{ej}

This lemma means that, as the polynomials $S_i$
are never identically zero, there is always a suitable
choice of the principal coefficients of the series involved
in a lift such that the tropicalization of this lift agrees
with the performed step of the tropical construction.
Moreover, the output of the step (namely the polynomials
$S_i$) can be considered as a single set $F_u$ for a later
construction. Clearly the input elements satisfy the
restrictions of the lemma, as their principal coefficients
are just different variables $c^j_i$. Thus, in the following
theorem, we use induction in order to show that a tropically
admissible construction agrees with the tropicalization of a
projective construction. The following is the main theorem of
the section.

\begin{thm}[General Lift]\label{Levantamientogeneral}
Suppose we are given the geometric construction of elements
$q_1,\ldots,q_s$ from elements $p_1,\ldots,p_n$. Suppose
that this construction can be meaningfully performed in the
projective space by Cramer's rule. If the construction of
each element $q_i$ is tropically admissible by Cramer's rule
then, for any specialization of the input tropical data given
by homogeneous coordinates $p_i = [p_i^1 : \ldots :
p_i^{m_i}]$, $1\leq i\leq n$, there exists a non empty set
$U$ in the space $(\mathbb{C}^*)^{m_1-1}\times \ldots\times
(\mathbb{C}^*)^{m_n-1}$ such that:
\begin{enumerate}
\item For every $(x_1,\ldots,x_n)$ in $U$ there exist
elements in the space of Puiseux series $P_1,\ldots,P_n$
such that $T(P_i)=p_i$, $Pc(P_i)=x_i$ and the
projective construction of $Q_1,\ldots,Q_s$ from
$P_1,\ldots,P_n$ is meaningful.
\item For all elements $P_1,\ldots,P_n$ in the
multiprojective space $\mathbb{P}^{m_1-1}(\mathbb{K})\times
\ldots\times \mathbb{P}^{m_n-1}(\mathbb{K})$ such that
$T(P_i)=p_i$ and $(Pc(P_1),$ $\ldots, Pc(P_n))$ $\in U$,
the tropicalization of the final elements of the construction
$Q_1,\ldots,Q_s$ agree with the tropical elements
$q_1,\ldots,q_s$ constructed using tropical determinants.
That is, all lifts with principal coefficients in $U$ yield
to the same tropical final elements.
\end{enumerate}
\end{thm}
\begin{proof}
We take generic projective lifts $P_1,\ldots,P_n$ with
$T(P_i)=p_i$, writing $Pc(P_i)=[c_i^1:\ldots:c_i^{m_i}]$,
indeterminate variables. Since each step of the construction
is given by Cramer's rule and all the variables are
different, we are in the hypotheses of
\ref{Levantamientounpaso}, taking for the first step
$f_i^j=c_i^j$, $F_i=\{f_i^1, \ldots ,f_i^{m_i}\}$,
$C_i=\{c_i^1,\ldots,c_i^{m_i}\}$. Each element of the
construction is admissible by Cramer's rule. Any intermediate
or output element is constructed from different objects. As
this element is tropically admissible, the condition of its
construction graph being a tree corresponds to the fact that
the input elements and hence the variables its parents depend
are different, so we are still in the conditions of
\ref{Levantamientounpaso}. This allows us to use induction in
\ref{Levantamientounpaso} because we will always have
disjoint sets of variables on the rows of our matrices. So,
all principal coefficients of all the steps in the
construction will be non-zero multihomogeneous polynomials in
the sets $C_i$.

We define $U$ as the subset of $(\mathbb{C}^*)^{m_1-1}\times
\ldots\times (\mathbb{C}^*)^{m_n-1}$ where all these
multihomogeneous polynomials do not vanish (considering
$(\mathbb{C}^*)^{m_i-1}\subseteq\mathbb{P}^{m_i-1}
(\mathbb{C})$ and taking homogeneous coordinates). If the
principal coefficients of the $P_i$ are in $U$, then we
obtain along the construction that all the principal
coefficients of the intermediate elements are non-zero. Then,
by lemma \ref{pseudocramerlema}, the tropicalization of each
step will be exactly the corresponding tropical determinant,
which is independent of the chosen lift $P_i$.

Of course, for $x_1,\ldots,x_n$ in $U$, one possible lift is
$P_i=[x_i^1t^{-p_i^1}:\ldots:x_i^{m_i}t^{-p_i^{m_i}}]$.
\end{proof}
\begin{defn}
Given a tropical geometric construction and a specialization
of the input data, we call \emph{general lift} of the input
data any lift whose principal coefficients belong to the set
$U$ defined above.
\end{defn}

\begin{rem}Theorem \ref{Levantamientogeneral} asserts that,
for every tropical geometric construction and for every input
data, if the construction graph of every element is a tree,
then there exists a lift that agrees with our tropical
construction, no matter what the input data is. Also, it is
remarkable that this theorem is stated in general dimension,
not just in the plane. So the result is valid for
constructions in $\mathbb{T}^n$, the only restriction we have
to consider is that of constructions involving only the stable
intersection of $n$ hyperplanes and the stable join of $n$
points.
\end{rem}

The following example shows what may happen if the
construction graph is not a tree and we are not in the
situation of theorem \ref{Levantamientogeneral}
\begin{ej}\label{ejemplociclo}
Suppose we are given $a,b,c$ three points in the plane. Let
$l_1=\overline{ab}$, $l_2=\overline{ac}$ be the lines
through these points and $p=l_1\cap l_2$. The construction of
$l_1$ and $l_2$ is tropically admissible by Cramer's rule,
but not the construction of $p$, because we have the cycle
$p,l_1,a,l_2,p$. The problem is that we have used twice the
point $a$ in order to construct $p$. Firstly it is used in
the construction of $l_1$ and then in the construction of
$l_2$.

So, after specialization, we may have some algebraic
relations making a pseudo-determinant identically zero for
every lift. For example, we take $a=[0:0:0]$, $b=[-2:1:0]$,
$c=[-1:3:0]$. Tropically, the construction yields
$l_1=[1:0:1]=1\odot x\oplus 0\odot y\oplus 1\odot z$,
$l_2=[3:0:3]=3\odot x\oplus 0\odot y\oplus 3\odot z$ and
finally $p=[3:4:3]=[0:1:0]\neq a$. But, for every lift of $a$,
$b$, $c$ such that the construction is well defined, the
final element must be the lift of $a$. These lifts take the
form $\widetilde{a}=[a_1:a_2:a_3]$, $\widetilde{b} = [b_1t^2 :
b_2t^{-1}:b_3]$, $\widetilde{c}=[c_1t:c_2t^{-3}:c_3]$, where
terms of bigger degree in the series do not affect the result.
In this case, $\widetilde{l}_1 = [-a_3b_2t^{-1}+a_2b_3 :
-a_1b_3+a_3b_1t2 : a_1b_2t^{-1}-a_2b_1t2] $ and
$\widetilde{l}_2 = [-a_3c_2t^{-3}+a_2c_3 :
-a_1c_3+a_3c_1t:a_1c_2t^{-3}+a_2c_1t]$ which tropicalize
correctly to $l_1$ and $l_2$ (as expected, because the
construction graphs of $l_1$ and $l_2$ are trees). Now, we
want to construct $\widetilde{p}$. Here, $O=\begin{pmatrix}
1&0&1\\2&0&2\end{pmatrix}$ and $A=\begin{pmatrix} -a_3b_2 &
-a_1b_3 & a_1b_2 \\-a_3c_2 & -a_1c_3 & a_1c_2
\end{pmatrix}$. Now it is easy to see that
$\Delta_{O^{2}} A^{2} = -a_1a_3b_2c_2+a_1a_3b_2c_2=0$. In
fact, $\widetilde{p}$ must be $\widetilde{a}$.

We observe that for $b$ and $c$ as in the example (which are
points in general position) the same lifting problem appears
for all $a=[r:s:0]$ with $r>-1$ and $s<1$. So the above
example is not at all an isolated case and it cannot be
avoided by perturbations of $a$, $b$ and $c$. These bad
conditioned cases arrive frequently when we are working with
non trivial constructions. So if we want to chain several
simple constructions we have to take these cases into
account, requiring some conditions on the construction graph
as formulated in the hypothesis of our theorem.
\begin{figure}
\begin{center}
\begin{picture}(70,80)
\put(-1,44){$a$}
\put(0,20){$b$}
\put(0,70){$c$}
\put(32,30){$l_1$}
\put(32,55){$l_2$}
\put(65,42){$p$}
\put(5,20){\line(2,1){25}}
\put(5,45){\line(2,-1){25}}
\put(5,45){\line(2,1){25}}
\put(5,70){\line(2,-1){25}}
\put(40,32){\line(2,1){25}}
\put(40,57){\line(2,-1){25}}
\end{picture}
\caption{The construction graph of $p$}\label{constp}
\end{center}
\end{figure}
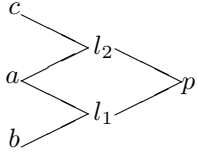
\end{ej}

\section{Constructive Pappus' Theorem}
Using theorem \ref{Levantamientogeneral}, we will now proof
the validity of the constructive version of Pappus' Theorem
proposed in \cite{RGST}

Let us start with the following specific lemma.

\begin{lem}\label{pappuslema}
If $k$ projective lines over the field of Puiseux series
\[L_i\equiv l_ix+l_i'y+l_i''z=0,\ l_il_i'l_i''\neq 0,
\ 1\leq i\leq k\] have a common point in the projective
plane, then the tropicalizations of all the lines contain a
common point in the tropical plane.
\end{lem}
\begin{proof}
The intersection point in the projective plane cannot be one
of $[1:0:0],\ [0:1:0],\ [0:0:1]$, because in that case one of
the coefficients of the equations would be zero, contrary to
the hypothesis. If the intersection point is in
$(\mathbb{K}^*)^2$, it is clear that the tropicalization of
all the lines contains the tropicalization of this point. It
remains to look what happens if the point is in one of the
coordinate lines. Suppose w.l.o.g. that the point is
$[0:a:b]$, $a\neq 0\neq b$. If we take the affine plane
corresponding to $\{z=1\}$, then the affine coordinates of
this point are $(0,a/b)$ and the affine equations of the
lines are of the form $L_i\equiv x=p_i(y-a/b)$, $p_i\neq 0$.
Let's take an element $x_0\neq 0$ such that
$o(x_0)>o(a/b)+\max\{o(p_i)|\ 1\leq i\leq k\}$.
There exists $y_i$ with $(x_0,y_i)\in L_i$.
$y_i=a/b+p_i^{-1}x_0$. As $o(p_i^{-1}x_0) =
o(x_0)-o(p_i)>o(a/b)$ we obtain that $o(y_i)=o(a/b),\ 1\leq
i\leq k$. So, the point $[-o(x_0):-o(a/b):0]\in T(L_i),\
1\leq i\leq k$. In fact, what we obtain is an infinite
number of points in the intersection of the tropical
lines.
\end{proof}

\begin{figure}
\begin{center}
\begin{picture}(70,80)
\put(0,10){1}
\put(0,25){2}
\put(0,40){4}
\put(0,55){3}
\put(0,70){5}
\put(32,32){$b$}
\put(32,62){$c'$}
\put(64,47){6}
\put(5,72){\line(4,-1){25}}
\put(5,57){\line(4,1){25}}
\put(5,42){\line(4,-1){25}}
\put(5,27){\line(4,1){25}}
\put(37,35){\line(2,1){25}}
\put(37,63){\line(2,-1){25}}
\put(69,47){\line(3,-1){35}}
\put(5,10){\line(4,1){100}}
\put(108,33){$a''$}
\end{picture}
\caption{The construction graph of $a''$}\label{consta2}
\end{center}
\end{figure}
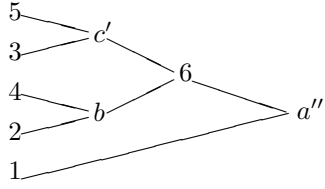

\begin{thm}[Constructive Pappus' Theorem]\label{Pappus}
Let 1, 2, 3, 4, 5 be five freely chosen points in the
tropical plane given by homogeneous coordinates. Define the
following additional three points and nine lines by a
sequence of stable join of points and stable meet of lines
operations (carried out by cross-products):
\[\begin{matrix}a=1\otimes 4,& b=2\otimes 4,& c=3\otimes 4,&
a'=1\otimes 5,& b'=2\otimes 5,& c'=3\otimes 5,\\6=b\otimes
c',& 7=a'\otimes c,& 8=a\otimes b',& a''=1\otimes 6,&
b''=2\otimes 7,& c''=3\otimes 8.\end{matrix}\] Then the
three tropical lines $a''$, $b''$, $c''$ are concurrent.
\end{thm}
\begin{proof}
Our goal is to see that the three constructed lines share a
common point. For each line, its construction is as follows
\[\begin{matrix}
a''=1\otimes((2\otimes 4)\otimes(3\otimes 5))\\
b''=2\otimes((3\otimes 4)\otimes(1\otimes 5))\\
c''=3\otimes((1\otimes 4)\otimes(2\otimes 5))\end{matrix}\]
We check on figure \ref{consta2} that the construction
graph of $a''$ is in fact a tree and the same holds for
$b''$ and $c''$. So for all general lifts of 1, 2, 3, 4, 5,
the construction of the lines $\widetilde{a}''$,
$\widetilde{b}''$, $\widetilde{c}''$ are well defined and
yields to the tropical lines $a''$, $b''$, $c''$. We notice
that we can not make in our tropical construction the stable
intersection of two of these lines, as they share input
points and cycles appear in the construction graph.
Nevertheless, the construction made in the projective space
satisfies the hypotheses of Pappus' theorem, so the lifts
$\widetilde{a}''$, $\widetilde{b}''$, $\widetilde{c}''$
must intersect in a common point in
$\mathbb{P}^2(\mathbb{K})$. Now, by \ref{pappuslema}, the
three tropical lines $a'',b'',c''$ must intersect.
\end{proof}

\begin{rem}
We have proved that our tropical construction is well
defined, but also that it agrees with the tropicalization
of almost all lifts of the construction. The lifts whose
principal coordinates are not in $U$ in theorem
\ref{Levantamientogeneral} include projective constructions
whose tropicalization is different to the one computed using
tropical Cramer's rule or non well defined constructions in
the projective space, such as the line passing through the
points $a$ and $a$. This later case does not appear in the
tropical situation if we interpret correctly the construction.

Suppose for example the case of Pappus' constructive version.
We start from five input points. Suppose that point 1 and 4
are the same. The line passing through 1 and 1 is not
tropically admissible, as there are cycles in the
construction graph. But the line passing through 1 and 4
is, even if 1=4, because we will take two different lifts
$\widetilde{1}\neq\widetilde{4}$ with
$T(\widetilde{1})=T(\widetilde{4})=1=4$. The tropical
construction is well defined for the whole space of
configurations of the original points, there is no need of
generality in the tropical space. The input elements of an
admissible tropical construction are completely free in the
sense that there exists no condition on these elements in
order to develop our construction and achieve the results.
\end{rem}

\section{Conclusions}
In this paper we have explored the possibility of developing
a tropical counterpart of classical geometric constructions.
In view of theorem \ref{Levantamientogeneral} we have
succeeded for Pappus' theorem. But we have also shown through
examples that there are several restrictions on the
constructions to apply this theorem. It would be interesting
to have some other remarkable examples of correct
tropicalization of classical theorems.

On the other hand, proposition \ref{consunpaso} shows that,
via Cramer's rule, the situation is specially simple for one
step constructions. We observe that this behavior is also
present in some other successful applications of tropical
geometry such as that of computing genus zero curves through
a general configuration of points developed in \cite{MIK2}.
This is, too, a ``one step mathematics'', since given a set
of tropically general points, we ``merely'' construct the
zero genus curves of given degree that passes through
these points.
The key problem seems to be handling tropical varieties
constructed from other previously constructed varieties.

So the morale suggested by the results presented in this
paper is that using tropical geometry is affordable (at this
moment) when dealing with ``one step mathematics'', but it
is not yet when dealing with geometric objects defined from
other objects that are not free in some algebraic sense.

\section*{Acknowledgments}
The author wants to thank Michel Coste, Ilia Itenberg and
Tomas Recio for very useful discussion, suggestions and
corrections.

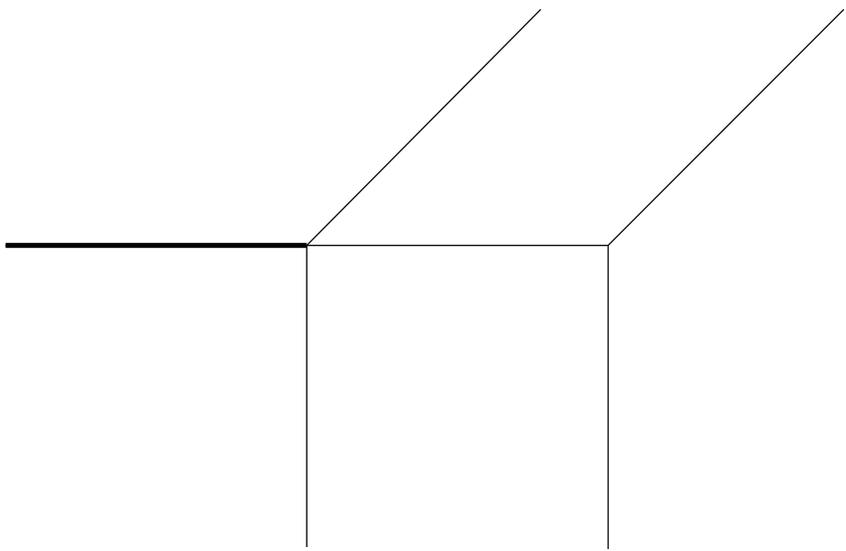
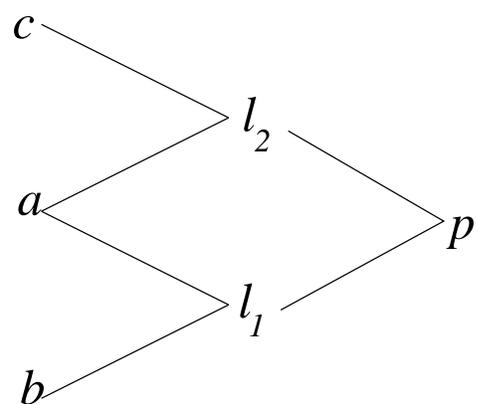
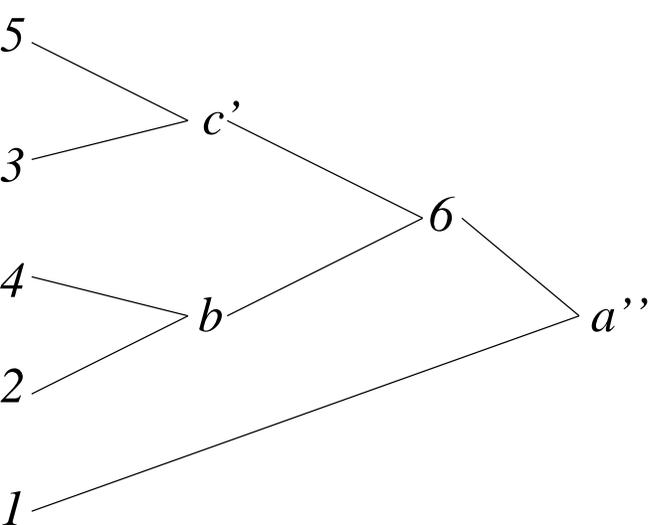
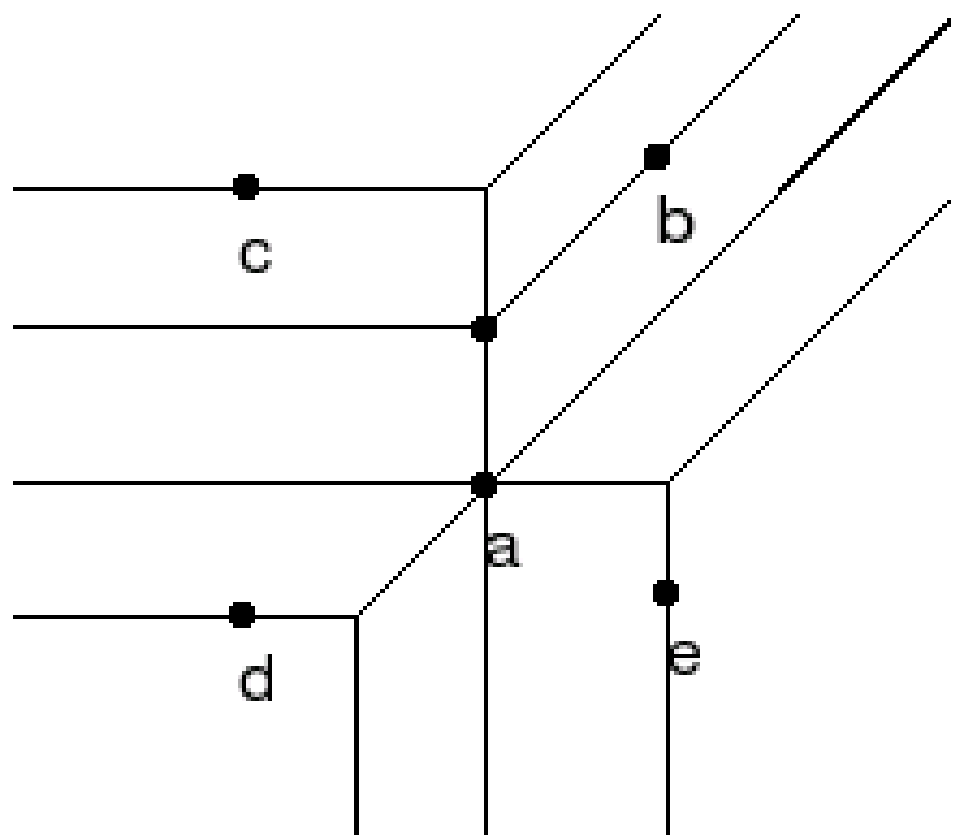
\end{document}